\documentclass{amsart}
\usepackage{amsmath, amssymb, amsthm}
\usepackage{url}
\usepackage[hidelinks]{hyperref}
\usepackage{mathtools}
\usepackage{upgreek}
\usepackage{environ}
\numberwithin{equation}{section}
\theoremstyle{plain}
\newtheorem{theorem}{Theorem}[section]
\newtheorem{lemma}[theorem]{Lemma}

\newtheorem*{definition*}{Definition}

\begin{document}
	\title[The Lattice Point Discrepancy of a Body of Revolution
	]{Omega Estimate for the Lattice Point Discrepancy of a Body of Revolution Using The Resonance Method} 
	
	\author{Nilmoni Karak}
	\address{Nilmoni Karak\\ Department of Mathematics \\
		Indian Institute of Technology Kharagpur\\
		Kharagpur-721302,  India.}
	\email{nilmonikarak@gmail.com,  nilmonimath@kgpian.iitkgp.ac.in}
	
	
	\maketitle
	
	\begin{abstract}
		Using a recent method developed by Mahatab, we obtain an improved $\Omega$-bound for the error term arising in lattice counting problem of bodies of revolution in $\mathbb R^3$ around a coordinate axis and having smooth boundary with bounded nonzero curvature. This strengthens an earlier result by K\"uhleitner and Nowak.
	\end{abstract}

\section{Introduction}
The classical Gauss circle problem is to find the number of points with integer coordinates inside a circle. This is equivalent to the number of representations of an integer as the sum of two squares. This problem extends to any closed domain in Euclidean spaces and has connections with number theory. An extensive study can be found in \cite{Krat}. As a special case, we may consider large homothetic smooth convex bodies.
%In lattice point theory, estimating the discrepancy between the volume and the number of lattice points in large homothetic smooth convex bodies constitutes a central problem.
Throughout this article we suppose $\mathfrak{B}\subset \mathbb{R}^3$ is a compact convex body with the origin $(0,0,0)\in \mathrm{int}(\mathfrak{B})$ and the boundary $\partial \mathfrak{B}$ is in class $C^{\infty}$ with positive Gaussian curvature. We consider the number of lattice points in the ``blown up" body $\sqrt{t}\mathfrak{B}$, for a large real $t$. Therefore, we define 
\begin{equation*}
	\mathcal{N}_\mathfrak{B}:= \#\{\mathbf{m}\in\mathbb{Z}^3: \mathbf{m}\in \sqrt{t}\mathfrak{B} \}.
\end{equation*}
Further, we define its \emph{lattice point discrepancy} as
\begin{equation}
	P_{\mathfrak{B}}(t) := 	\mathcal{N}_\mathfrak{B} - \text{vol}(\mathfrak{B})t^{3/2}.
\end{equation}
In this direction, Hlawka \cite{Hlawka} was the first to study the asymptotic behavior of $\mathcal{N}_{\mathfrak{B}}$ for a general convex body $\mathfrak{B}$. Assuming the body $\mathfrak{B}$ to be a body of revolution with respect to any one of the coordinate axes, Chamizo \cite{Cham} showed that
\begin{align*}
	P_{\mathfrak{B}}(t)\ll t^{11/16}.
\end{align*}
Extensions of the above asymptotic result to higher dimensions can be found in \cite{Gu,Mul98}. In the same paper, Hlawka \cite{Hlawka} also studied $\Omega$-bound for $P_{\mathfrak{B}}$. Through a series of papers, K\"uhleitner and Nowak \cite{K00, KN, Nowak85}\footnote{Throughout the article, we use $\log_rx=\log\log\dots\log x\,\,\,(r\, \text{times})$.
} improved the $\Omega$ bound. The final improvement in \cite{KN} is obtained using Soundararajan's \cite[Lemma 2.1]{Sound} method, which proves
\begin{equation}\label{K-Nowak}
		P_{\mathfrak{B}}(t)= \Omega_{-}\Big(t^{1/2}(\log t)^{1/3}(\log_2 t)^{\frac{2}{3}(\sqrt{2}-1)}(\log_3 t)^{-2/3}\Big).
\end{equation}
Recently, Mahatab \cite{K. Mahatab} has obtained a sharper $\Omega$-bound for the circle problem and the divisor problems using the resonance method, inspired from a resonator used by Aistleitner, Mahatab and  Munsch \cite{AMM}.  In this paper, we use Mahatab's resonator \cite[Theorem 6]{K. Mahatab} to obtain an improved lower bound for \eqref{K-Nowak}. In fact, we improve on the power of $\log_3 x$ from $-2/3$ to $-1/3$.
\begin{theorem}\label{Main result}
Let $\mathfrak{B}$ be a compact, convex body in $\mathbb{R}^3$ whose interior contains the origin and is invariant under rotations around one of the coordinate axes. Assume that its boundary $\partial{\mathfrak{B}}$ is  in $C^{\infty}$ with positive bounded curvature.	Then for a large $t$, we have
	\begin{align*}
		P_{\mathfrak{B}}(t)= \Omega_{-}\Big(t^{1/2}(\log t)^{1/3}(\log_2 t)^{\frac{2}{3}(\sqrt{2}-1)}(\log_3 t)^{-1/3}\Big).
	\end{align*}
\end{theorem}
In the next section, we collect all required tools to prove Theorem \ref{Main result}.

\section{Preliminaries}
	First, we mention a theorem, due to Mahatab \cite{K. Mahatab}, as the most important tool for our result. Let $(\lambda_n)_{n=1}^{\infty}$ be a non-decreasing sequence of positive real numbers  and $\alpha$ be a positive real parameter. We consider a finite linearly independent set $\mathcal{M}$ over $\mathbb{Q}$ such that $\mathcal{M}\subseteq \{\lambda_n: C_0 \alpha \leq \lambda_n\leq 2\alpha\}$, with $0<C_0<2$.  We denote the cardinality of the set $\mathcal{M}$ by $M$  and fix four real numbers satisfying $0<A_4<A_3<A_2<A_1$. 
	%With all the notations in our hand, we have the following result due to Mahatab \cite{K. Mahatab}:
	\begin{theorem}[K. Mahatab] \label{K-M Theorem}
		Let $(a_n)_{n=1}^{\infty}$ be a sequence of positive real numbers. Then for large $T$, we have
		\begin{align*}
			\max_{T^{A_3/2}\leq t\leq 2A_2^2T^{A_2}\log^2 T}\left( \sum_{n\leq T^{A_1}} a_n \cos(t\lambda_n)\right)\geq \frac{\pi}{4e}\sum_{n\in \mathcal{M}}a_n &+ O\left( T^{A_3-A_2}e^{2M/C_0}\left(\sum_{\lambda_n\leq 4\alpha}a_n\right)\right) \\
			& + O\left(\frac{T^{-A_4}}{\alpha}\sum_{n\leq T^{A_1}}a_n\right).
		\end{align*}
	\end{theorem}
Now we define the \emph{tac function} $H$ for a convex body $\mathfrak{B}$ satisfying the conditions defined before.
\begin{definition*} \label{Tac function}
	For $\mathbf{u}\in\mathbb{R}^3$, the tac function $H$ of a convex body $\mathfrak{B}$ is defined by 
	\begin{align*}
		H(\mathbf{u})=\max_{\mathbf{u}\in\mathfrak{B}}\langle \mathbf{u},\mathbf{v}\rangle,
	\end{align*}
	where $\langle , \rangle$ denotes the usual Euclidean inner product. The following properties are immediate.
	\begin{enumerate}
		\item $H$ is a real positive homogeneous function of degree 1.
		\item For any $\mathbf{u}\in \mathbb{R}^3$, there exist positive constants $a$ and $b$ such that
		\begin{align}\label{tac prprty 1}
			a\,\lVert \mathbf{u}\rVert\leq H(\mathbf{u})\leq b\, \lVert \mathbf{u}\rVert,
		\end{align} 
		where $\lVert\cdot\rVert$ denotes the Euclidean norm.
		\item For $(u_1,u_2,u_3)\in\mathbb{R}^3$,
		\begin{align}\label{tac prprty 2}
			H(u_1,u_2,u_3)=H\left(\sqrt{u_1^2+u_2^2},0,u_3\right).
		\end{align}
	\end{enumerate}
\end{definition*}		
The proof of Theorem \ref{Main result} also uses the next lemma, which can be obtained by taking $s=3$ in \cite[Eq. (13)]{Nowak-1}. For more information, the reader is suggested to see \cite[Section 2]{KN}. 
%In the next lemma, we collect another important tool for the proof of our main result.
	\begin{lemma}\label{Borel lemma}
		For a body $\mathfrak{B}$ satisfying the conditions stated before and a large real parameter $t$, the Borel mean value of the lattice point discrepancy $P_{\mathfrak{B}}$ is defined as
		\begin{equation}\label{Borel}
			B(t):= \frac{1}{\Gamma(k+1)}	\int_{0}^{\infty}e^{-u}u^k P_{\mathfrak{B}}\big(Xu\big) \mathrm{d}u,
		\end{equation}
		where $X=X(t)= 1/(\log t)$ and $k=k(t)= t^2\log t$. Then $B(t)$ has the following asymptotic formula:
		\begin{equation} \label{Formula of B(t)}
			B(t) =- \frac{1}{2\pi} tS(t) +O(t^{3/4+\epsilon}),
		\end{equation}
		with
		\begin{equation}\label{S(t)}
				S(t) := \sum_{0<\lVert {\bf m}\rVert\leq t^{\epsilon}\sqrt{\log t}} \frac{\theta({\bf m})}{\lVert {\bf m}\rVert^2} \exp\Big(-\frac{1}{2}\pi^2 X H\big({\bf m}\big)^2\Big)\cos\Big(2\pi H\big({\bf m}\big)t\Big),
		\end{equation}
	where the coefficients $\theta({\bf m})$ are positive and uniformly bounded from above and ${\bf m}\in\mathbb{Z}^3$.	
	\end{lemma}
	
\section{Proof of Theorem 1.1}
For ${\bf m}=(m_1,m_2,m_3)\in\mathbb{Z}^3$, we set $l=m_1^2+m_2^2$ and construct a bijection $\uppsi:\mathbb{N} \rightarrow\mathbb{N}\times \mathbb{Z}$ as $n\mapsto\uppsi(n)= (l,m_3)$, such that the sequence $(\lambda_n)_{n=1}^{\infty}$ defined by
\begin{equation*}
	\lambda_n:=H({\bf m})\Big|_{(l,m_3)=\uppsi(n)}= H(\sqrt{l},0,m_3)\Big|_{(l,m_3)=\uppsi(n)}
\end{equation*} is non-decreasing.
Implementing \eqref{tac prprty 2} in \eqref{S(t)}, one has
		\begin{equation*}
				S(t) := \sum_{0<l+m_3^2\leq t^{2\epsilon}\log t} \frac{f(l,m_3)}{l+m_3^2} \exp\Big(-\frac{1}{2}\pi^2 X H\big(\sqrt{l},0,m_3\big)^2\Big)\cos\Big(2\pi H\big(\sqrt{l},0,m_3\big)t\Big),
			\end{equation*}
where
\begin{align}\label{r(l)}
	f(l,m_3):=\sum_{m_1^2+m_2^2=l} \theta (m_1,m_2,m_3)\asymp r(l),
\end{align}			
and $r(l)$ denotes the number of ways of writing $l$ as sum of two squares. Now we define a sequence $\left(a_n\right)_{n=1}^{\infty}$ as
\begin{equation}\label{a_n}
	a_n := 
	\begin{cases}
		\displaystyle\frac{f(l, m_3)}{l + m_3^2} \exp\left(-\frac{1}{2}\pi^2 X H\big(\sqrt{l}, 0, m_3\big)^2\right)\bigg|_{(l, m_3) = \uppsi(n)}, & \text{if } l + m_3^2 \leq t^{2\epsilon}\log t, \\
		0, & \text{otherwise}.
	\end{cases}
\end{equation}

We assume a set $\mathcal{M}$ such that $\{ \lambda_n: n\in \mathcal{M}\} \subset \left[C_0 \alpha, 2\alpha\right]$, where $0<C_0<2$. Also, we choose $A_1=2$, $A_2=3/2, A_3=1$ and $A_4=7/8$. By Theorem \ref{K-M Theorem}, there exists $t$ satisfying $\sqrt{T}\leq t\leq 5T^{3/2}\log^2 T$ such that 
\begin{equation} \label{lower bound of S(t)}
	S(t)\geq \cfrac{\pi}{4e} \sum_{n\in \mathcal{M}}a_n + O\left( T^{-1/2}e^{2M/C_0}\left(\sum_{\lambda_n\leq 4\alpha}a_n\right)\right) +  O\left(\frac{T^{-7/8}}{\alpha}\sum_{n\leq T^{2}}a_n\right),
\end{equation}
where $\alpha$ is a real parameter to be determined later.

Since the \emph{tac function} $H$ is homogeneous, there exist constants $0<c_1<c_2$ and $0<c_3<c _4$ depending on $\mathfrak{B}$ such that the interval $[c_1,c_2]\times [c_3,c_4]$ in $(u_1,u_3)$-plane lies between the curves $H(u_1,0,u_3)=C_0$ and $H(u_1,0,u_3)=2$. Therefore, the condition $(\sqrt{l},m_3) \in [c_1\alpha,c_2\alpha]\times [c_3\alpha,c_4\alpha]$ implies $H(\sqrt{l},0,m_3)\in [C_0 \alpha, 2\alpha]$ for integers $l>0$ and $m_3$.

 Now we construct the following resonating set $\widehat{\mathcal{M}}$ such that $\mathcal{M}$ is the preimage of $\widehat{\mathcal{M}}$ under the map $\uppsi$:
%\begin{equation*}
$$\widehat{\mathcal{M}} = \left\{
%\begin{aligned}
(l, m_3) \in \mathbb{N}^2 :\  c_1^2 \alpha^2 \leq l \leq c_2^2 \alpha^2, \ c_3 \alpha \leq m_3 \leq c_4 \alpha, \ \omega(l)=\left[ \lambda \log_2 \alpha \right] \ \text{and}\ l\in \mathcal{A}
%\end{aligned}
\right\},$$
%\end{equation*}
where $\mathcal{A}:=\left\{q\in\mathbb{N}: p\equiv 1 (\hspace{-1.9 mm}\mod4), \text{if}\ \,p|q\,\ \text{and}\ p\ \text{is prime};  \text{and}\ q\ \text{is square free}\right\}$, $\omega(l)$ denotes the number of distinct prime factors of $l$ and $\lambda$ is a positive constant. 

Now using \eqref{r(l)} and the assumption $XH(\sqrt{l}, 0, m_3)^2\ll 1$ in \eqref{a_n}, we have
\begin{align}\label{sum a_n}
\nonumber	\sum_{n\in \mathcal{M}}a_n &\gg \cfrac{1}{\alpha^2}\sum_{c_3\alpha\leq m_3\leq c_4\alpha} \sum_{\substack{c_1^2\alpha^2\leq l\leq c_2^2\alpha^2\\
	\omega(l)= [\lambda\log_2 \alpha],\ l\in \mathcal{A}}} r(l)\\
&\gg \cfrac{1}{\alpha} \sum_{\substack{c_1^2\alpha^2\leq l\leq c_2^2\alpha^2\\
		\omega(l)= [\lambda\log_2 \alpha],\ l\in \mathcal{A}}} r(l).
\end{align}
Let $S$ be the cardinality of the set
\begin{align*}
	S_{\alpha, \Lambda}= \{l\in \mathbb{N}: a_1^2\alpha\leq l\leq a_2^2\alpha^2, l\in \mathcal{A}\ \text{and}\ \omega(l)=\Lambda \}. 
\end{align*}
Using Sathe's \cite{Sathe} (See also \cite[Section II.6]{Tenenbaum}) and Stirling's formula, one can obtain 
\begin{equation*}
	S\asymp \frac{\alpha^2}{\log \alpha} \frac{(\frac{1}{2}\log_2 \alpha)^{\Lambda -1}}{(\Lambda-1)!} \asymp  \frac{\alpha^2}{\sqrt{\log_2 \alpha}}(\log \alpha)^{\lambda -1 -\lambda\log\lambda - \lambda \log 2},
\end{equation*}
where we have chosen $\Lambda = [\lambda\log_2\alpha]$ and therefore
\begin{equation*}
	M=|\mathcal{M}| =|\mathcal{\widehat{M}}| \asymp \frac{\alpha^3}{\sqrt{\log_2 \alpha}}(\log \alpha)^{\lambda -1 -\lambda\log\lambda - \lambda \log 2}.
\end{equation*}
Note that, for $l\in \mathcal{A}$, $r(l)\geq 2^{w(l)} = 2^{[\lambda \log_2 \alpha]}\gg (\log \alpha )^{\lambda \log 2}$. Using this lower bound of $r(l)$ in \eqref{sum a_n}, we have
\begin{align}\label{sum a_n 2}
	\nonumber\sum_{n\in \mathcal{M}}a_n&\gg \frac{(\log \alpha)^{\lambda \log 2}}{\alpha} \frac{\alpha^2}{\sqrt{\log_2 \alpha}}(\log \alpha)^{\lambda -1-\lambda \log \lambda -\lambda \log 2}\\
	& = \frac{\alpha}{\sqrt{\log_2 \alpha}}(\log \alpha)^{\lambda -1-\lambda \log \lambda}.
\end{align}
Now we choose $\alpha$ such that $M$ is of order $\log T$. Indeed, 
\begin{equation*}
	\alpha = \frac{1}{C}(\log T)^{1/3}(\log_2 T)^{\frac{1}{3}(1-\lambda +\lambda\log\lambda +\lambda \log2)}(\log_3 T)^{1/6},
\end{equation*}
where C is a large positive constant. Note that $\log\alpha \asymp \log_2 T$ and $\log_2 \alpha \asymp \log_3 T$.   Since $X\ll (\log T)^{-1}$ and $H(\sqrt{l},0,m_3)\ll \alpha$, the assumption $XH(\sqrt{l}, 0, m_3)^2\ll 1$ is now verified.
Substituting the above choice of $\alpha$ in \eqref{sum a_n 2},
\begin{align*}
	\sum_{n\in \mathcal{M}}a_n \gg (\log T)^{1/3} (\log_2 T)^{\frac{2}{3}(\lambda-1-\lambda\log \lambda)+\frac{1}{3}\lambda \log 2}(\log_3 T)^{-1/3}.
\end{align*}
To optimize the exponent of $\log_2 T$, we choose $\lambda=\sqrt{2}$ and we obtain
\begin{align}\label{Final lower bound of sum a_n}
	\sum_{n\in \mathcal{M}}a_n \gg (\log T)^{1/3} (\log_2 T)^{\frac{2}{3}(\sqrt{2}-1)}(\log_3 T)^{-1/3}.
\end{align}
To estimate the error, first we note that 
\begin{align*}
	M\asymp \frac{\log T}{C^3}.
\end{align*}
We choose $C$ to be large enough satisfying 
\begin{align*}
	e^{2M/C_0}\ll T^{1/4-\epsilon}.
\end{align*}
Now, using \eqref{tac prprty 1} and \eqref{tac prprty 2}, we calculate\footnote{The quantity $r_3(n)$ denotes the number of ways of writing $n$ as sum of three squares.}
\begin{align*}
	\sum_{\lambda_n\leq 4\alpha}a_n&\ll \sum_{0<H(\sqrt{l},0,m_3)\leq 4\alpha} \frac{r(l)}{l+m_3^2}\\
	&\leq \sum_{0<a\lVert {\bf m}\rVert\leq 4\alpha} \frac{1}{\lVert \bf{m}\rVert^2}\\
	&= \sum_{n\leq (16\alpha^2)/a^2} \frac{r_3(n)}{n}\\
	&= \int_{1}^{ (16\alpha^2)/a^2} \frac{1}{x}\ \mathrm{d}\left(\sum_{n\leq x}r_3(n)\right)\\
	& \ll\int_{1}^{ (16\alpha^2)/a^2} \frac{1}{x^{1/2}}\mathrm{d} x \ll \alpha \ll T^{\epsilon},
\end{align*}
where we use the bound $\sum_{n\leq x}r_3(n)\ll x^{3/2}$.

Similarly,
\begin{align*}
	\sum_{n\leq T^2}a_n &\ll \sum_{0<\lVert {\bf{m}} \rVert\leq t^{\epsilon}\sqrt{\log t} }\frac{1}{\lVert \mathbf{m} \rVert^2}\\
	&= \int_{1}^{t^{2\epsilon_0}\log t} \frac{1}{x}\ \mathrm{d}\left(\sum_{n\leq x}r_3(n)\right)\\
	&\ll t^{\epsilon}\sqrt{\log t}\ll T^{3\epsilon
	}.
\end{align*}
The above two estimates show that the two error terms in the right hand side of \eqref{lower bound of S(t)} are smaller than the main term \eqref{Final lower bound of sum a_n}. In fact,
\begin{align*}
	\text{Error}&\ll  T^{-1/2}e^{2M/C_0}\left(\sum_{\lambda_n\leq 4\alpha}a_n\right) +  \frac{T^{-7/8}}{\alpha}\left(\sum_{n\leq T^{2}}a_n\right)\\
	&\ll T^{-1/4+\epsilon}+ T^{-7/8+\epsilon}\ll T^{-1/4+\epsilon}.
\end{align*}
Using \eqref{Final lower bound of sum a_n} in \eqref{lower bound of S(t)}, it follows that for large $T$, there exists $t$ satisfying $t\in \sqrt{T}\leq t\leq 5T^{3/2}\log^2 T$ such that
\begin{align*}
	S(t)\gg (\log t)^{1/3} (\log_2 t)^{\frac{2}{3}(\sqrt{2}-1)}(\log_3 t)^{-1/3}.
\end{align*}
Then from \eqref{Formula of B(t)}, we have
\begin{align} \label{lower bound of -B(t)}
	-B(t)\gg t(\log t)^{1/3} (\log_2 t)^{\frac{2}{3}(\sqrt{2}-1)}(\log_3 t)^{-1/3}.
\end{align}
Let 
\begin{align*}
	\mathcal{F}(w):=(\log w)^{1/3} (\log_2 w)^{\frac{2}{3}(\sqrt{2}-1)}(\log_3 w)^{-1/3}.
\end{align*}
We will prove our $\Omega$ result by the method of contradiction.
Suppose for any  $\delta>0$, there exists a constant $C_1$ such that
\begin{align*}
	-P_{\mathfrak{B}}(w)\leq C_1 +\delta w^{1/2}\mathcal{F}(w),
\end{align*}
for all $w>0$.
Now applying \eqref{Borel} to get
\begin{align*}
	-B(t)&\leq \dfrac{1}{\Gamma(k+1)}\int_{0}^{\infty} e^{-w} w^k \left(C_1+\delta w^{1/2}\mathcal{F}(w)\right) \mathrm{d}w\\ &= C_1 + \frac{\delta}{\Gamma(k+1)}\int_{0}^{\infty} e^{-w}w^k(Xw)^{1/2}\mathcal{F}(w) \mathrm{d}w,
\end{align*}
for all $t>0$. Estimating the above infinite integral by Hafner's result \cite[Lemma 2.3.6]{Hafner}, we obtain that
\begin{align*}
	-B(t)\leq C_1 +C_2\delta (kX)^{1/2} \mathcal{F}(kX)\leq C_1+C_2\, \delta \, t\, \mathcal{F}(t^2),
\end{align*}
for some constant $C_2>0$ and for any $\delta>0$.
However, this contradicts \eqref{lower bound of -B(t)}, which proves our required $\Omega$ result.

\section*{Acknowledgment} The author would like to thank Dr.~Kamalakshya Mahatab for his valuable guidance and many helpful discussions throughout this work. The work of the author is supported by Prime Minister's Research Fellowship (PMRF) with PMRF Id: 2403449.

\end{document}